\documentclass[11pt,a4paper]{article}
\pdftrailerid{}

\usepackage{microtype}
\usepackage[utf8]{inputenc}
\usepackage{setspace}
\usepackage{fullpage}
\usepackage{amsmath}
\usepackage{amssymb}
\usepackage{amsthm}
\usepackage{verbatim}
\usepackage{enumerate}
\usepackage{amsfonts}
\usepackage{amsmath}
\usepackage[backref]{hyperref}
\hypersetup{%
    colorlinks,
    linkcolor={red!60!black},
    citecolor={green!60!black},
    urlcolor={blue!60!black},
}
\usepackage{comment}
\usepackage{cleveref}

\usepackage{enumitem}
\def\rmlabel{\upshape({\itshape \roman*\,})}

\makeatletter
\let\@@pmod\pmod
\DeclareRobustCommand{\pmod}{\@ifstar\@pmods\@@pmod}
\def\@pmods#1{\mkern4mu({\operator@font mod}\mkern 6mu#1)}
\makeatother

\usepackage[normalem]{ulem}

\usepackage{tikz}
\usetikzlibrary{calc,positioning,decorations.pathmorphing,decorations.pathreplacing}
\usepackage{ifthen}

\newtheorem{theorem}{Theorem}
\newtheorem{conj}[theorem]{Conjecture}
\newtheorem{lemma}[theorem]{Lemma}
\newtheorem{prop}[theorem]{Proposition}
\newtheorem{fact}[theorem]{Fact}

\newtheorem{cor}[theorem]{Corollary}
\newtheorem*{remark*}{Remark}
\newtheorem{claim}[]{Claim}

\theoremstyle{definition}

\newtheorem{definition} [theorem] {Definition}

\DeclareMathOperator{\Bin}{{\rm Bin}}
\let\Re\relax
\DeclareMathOperator{\Re}{{\rm Re}}

\def\[{\left(}
\def\]{\right)}

\def\llceil{\left\lceil}
\def\rrceil{\right\rceil}

\let\::%

\newcommand{\EE}{\mathbb{E}}
\newcommand{\PP}{\mathbb{P}}

\usepackage{dsfont}

\let\epsilon=\varepsilon

\newcommand{\quot}[2]{\mathchoice%
	{\left.\raisebox{.1em}{$\displaystyle{#1}$}\kern-1pt/\raisebox{-.2em}{$\displaystyle{#2}$}\right.}
	{\left.\raisebox{.1em}{${#1}$}\kern-1pt/\raisebox{-.2em}{${#2}$}\right.}
	{\left.\raisebox{.1em}{$\scriptstyle{#1}$}\kern-1pt/\raisebox{-.2em}{$\scriptstyle{#2}$}\right.}
{\left.\raisebox{.1em}{$\scriptscriptstyle{#1}$}\kern-1pt/\raisebox{-.2em}{$\scriptscriptstyle{#2}$}\right.}
}

\usepackage{datetime}

\usepackage[abbrev,msc-links,backrefs]{amsrefs}
\usepackage{doi}

\renewcommand{\PrintDOI}[1]{\doi{#1}}

\usepackage[mathlines]{lineno}
\usepackage{etoolbox} %

\newcommand*\linenomathpatch[1]{%
	\expandafter\pretocmd\csname #1\endcsname {\linenomath}{}{}%
	\expandafter\pretocmd\csname #1*\endcsname{\linenomath}{}{}%
	\expandafter\apptocmd\csname end#1\endcsname {\endlinenomath}{}{}%
	\expandafter\apptocmd\csname end#1*\endcsname{\endlinenomath}{}{}%
}
\newcommand*\linenomathpatchAMS[1]{%
	\expandafter\pretocmd\csname #1\endcsname {\linenomathAMS}{}{}%
	\expandafter\pretocmd\csname #1*\endcsname{\linenomathAMS}{}{}%
	\expandafter\apptocmd\csname end#1\endcsname {\endlinenomath}{}{}%
	\expandafter\apptocmd\csname end#1*\endcsname{\endlinenomath}{}{}%
}

\expandafter\ifx\linenomath\linenomathWithnumbers
\let\linenomathAMS\linenomathWithnumbers
\patchcmd\linenomathAMS{\advance\postdisplaypenalty\linenopenalty}{}{}{}
\else
\let\linenomathAMS\linenomathNonumbers
\fi

\linenomathpatchAMS{gather}
\linenomathpatchAMS{multline}
\linenomathpatchAMS{align}
\linenomathpatchAMS{alignat}
\linenomathpatchAMS{flalign}
\linenomathpatch{equation}

\title{The $\!{}\bmod k$ chromatic index of random graphs%
	\thanks{\scriptsize This research has been partially supported by
		Coordenação de Aperfeiçoamento  
		de Pessoal de Nível Superior - Brasil -- CAPES -- Finance Code 001. 
		F.~Botler is partially supported by CNPq (423395/2018-1)
		and by FAPERJ (211.305/2019).
		L.~Colucci is supported by FAPESP (2020/08252-2).
		Y.~Kohayakawa is partially supported by
		CNPq (311412/2018-1, 423833/2018-9, 406248/2021-4) and
                FAPESP (2018/04876-1, 2019/13364-7). 
		The research that led to this paper started at WoPOCA 2019, which was
		financed by FAPESP (2015/11937-9) and CNPq (425340/2016-3,
		423833/2018-9). FAPERJ is the Rio de Janeiro Research Foundation.
		FAPESP is the S\~ao Paulo Research Foundation.  CNPq is the National 
		Council for Scientific and Technological Development of Brazil.}}
\author{Fábio Botler,\footnote{\scriptsize Programa de Engenharia de
		Sistemas e Computação, COPPE, Universidade Federal do Rio de
		Janeiro, Brazil} 
	$\;$Lucas Colucci$^\ddag$ and Yoshiharu
        Kohayakawa\footnote{\scriptsize
		Instituto de Matemática e Estatística, Universidade de São Paulo,
		Brazil}}
\date{}

\begin{document}
\shortdate
\yyyymmdddate
\settimeformat{ampmtime}
\onehalfspace
\maketitle

\begin{abstract}
  The \emph{mod \(k\) chromatic index} of a graph \(G\) is the minimum
  number of colors needed to color the edges of $G$ in a way that the
  subgraph spanned by the edges of each color has all degrees
  congruent to~$1\pmod*k$.  Recently, the authors proved that the
  mod~\(k\) chromatic index of every graph is at most \(198k-101\),
  improving, for large~$k$, a result of Scott [{\em Discrete
    Math.~175}, 1--3 (1997), 289--291].  Here we study the mod~\(k\)
  chromatic index of random graphs.  We prove that for every integer
  $k\geq2$, there is $C_k>0$ such that if $p\geq C_kn^{-1}\log{n}$ and
  $n(1-p) \rightarrow\infty$ as $n\to\infty$, then the following
  holds: if~\(k\) is odd, then the mod~$k$ chromatic index of~$G(n,p)$
  is asymptotically almost surely equal to~$k$, while if~\(k\) is
  even, then the mod~\(k\) chromatic index of $G(2n,p)$
  (respectively~$G(2n+1,p)$) is asymptotically almost surely equal to
  $k$ (respectively~\(k+1\)).
\end{abstract}

\section{Introduction}

Throughout this paper, all graphs are simple and~$k\geq2$ is a fixed
integer.  If $G=(V,E)$ is a graph and~$F\subset E$, then the subgraph
of~$G$ \textit{spanned} by~$F$ is $G[F]=(W,F)$, where~$W$ is the set
of vertices of~$G$ that are incident to at least one edge in~$F$.
Note that, in particular, $G[F]$ has no isolated vertices.  A
\emph{\(\chi'_k\)-coloring} of a graph \(G\) is a coloring of the
edges of \(G\) in which the subgraph spanned by the edges of each
color has all degrees congruent to~$1\pmod*k$.  The \emph{mod \(k\)
  chromatic index} of~\(G\), denoted~$\chi'_k(G)$, is the minimum
number of colors in a \(\chi'_k\)-coloring of~\(G\).  Note that
isolated vertices play no role in these definitions, and hence
$\chi_k'(G)=\chi_k'(G-v)$ if~$v$ is an isolated vertex in~$G$.  Since
a proper coloring of the edges of~$G$ is a $\chi'_k$-coloring, we have
$\chi'_k(G)\leq\chi'(G)$, where, as usual, $\chi'(G)$~denotes the
chromatic index of~$G$.  It turns out that~$\chi_k'(G)$ is
usually much smaller than~$\chi'(G)$.  In~1991,
Pyber~\cite{pyber1991covering} proved that $\chi'_2(G) \leq 4$ for
every graph $G$, and in~1997 Scott~\cite{scott1997graph} proved that
$\chi'_k(G) \leq 5k^2\log k$ for every graph~\(G\).  Recently, the
authors~\cite{BoCoKo20} proved that $\chi'_k(G)\leq198k-101$ for
any~$G$,\footnote{By making use of a recent result of Hasanvand (see
  the comments following Theorem~1.4 in~\cite{hasanvand22:_modul}),
  the proof in~\cite{BoCoKo20} yields \(\chi'_k(G) \leq 38k -37\).}  which improves Scott's bound
for large~$k$ and is sharp up to the multiplicative constant, as every
graph with a vertex of degree~$k$ requires at least~$k$ colors in any
\(\chi'_k\)-coloring.  In this paper, we study the behavior of the
mod~\(k\) chromatic index of the random graph~$G(n,p)$ for a wide
range of~$p=p(n)$.  More specifically, we prove the following result.

\begin{theorem}\label{thm:main}
  For every integer $k \geq 2$, 
  there is a constant~$C_k$ such that if $p \geq C_kn^{-1}\log{n}$ and
  $n(1-p) \rightarrow\infty$, then the  following holds as $n \to \infty$.
  \begin{enumerate}[label=\rmlabel]
  \item If $k$ is even, then
    \begin{equation*}
      \lim_{n\to\infty}\PP\big(\chi'_k(G(2n,p)) = k\big)=1,
    \end{equation*}
    and
    \begin{equation*}
      \lim_{n\to\infty}\PP\big(\chi'_k(G(2n+1,p)) = k+1\big)=1.                    
    \end{equation*}
  \item If $k$ is odd, then
    \begin{equation*}
      \lim_{n\to\infty}\PP\big(\chi'_k(G(n,p)) = k\big)=1.
    \end{equation*}
  \end{enumerate}
\end{theorem}

Theorem~\ref{thm:main} extends a result of the authors
\cite{etc} that dealt with the case $k=2$ and
$C\sqrt{n^{-1}\log{n}} < p < 1-1/C$ for some $C$. 

In Section~\ref{sec:technical-lemmas}, we present some technical
lemmas that we shall need and in Section~\ref{sec:main-theorem} we
prove Theorem~\ref{thm:main}.  While Theorem~\ref{thm:main} tells us
that the typical value of the mod~$k$ chromatic index is at most~$k+1$
for a wide range of~$ p=p(n)$, in
Section~\ref{sec:deterministic-lower-bound}, we give a sequence
\(G_2,G_3,\ldots\) of graphs for which \(\chi'_k(G_k) \geq k+2\) for
every $k\geq2$ (see Proposition~\ref{prop:lwbd}).
Section~\ref{sec:concluding-remarks} contains some remarks on other
ranges of~$p$.

A property~$P$ holds for~$G(n,p)$ \textit{asymptotically almost
  surely} (a.a.s.) or \textit{with high probability} if the
probability that~$G(n,p)$ satisfies~$P$ tends to~$1$ as~$n\to\infty$.
The asymptotic notation~$o(1)$, $\ll$ and~$\gg$ will always be with
respect to~$n\to\infty$.

The main results of this work were announced in the extended
abstract~\cite{botler21:_mod_chrom_index_rgs}. 

\section{Technical lemmas}
\label{sec:technical-lemmas}
	
Given a graph $G$ and \(i\in \{1,\ldots,k\}\), let \(V_i=V_i(G)\) be
the set of vertices~$v$ of~\(G\) with $\deg_G(v)\equiv i\pmod*k$.
Let~$n_i(G)=|V_i|$ and let $G_i=G[V_i]$.  The sets \(V_1,\ldots, V_k\)
are the \emph{degree classes} of~\(G\).

In this section, after presenting some auxiliary probabilistic and
random graph results, we prove that~$G=G(n,p)$ in a wide range
of~\(p\) with high probability satisfies (\textit{a})~$n_i(G)$ do not
deviate much from~$n/k$ (Lemma~\ref{lem:bal_deg_cls}),
(\textit{b})~every vertex of~$G$ has about~$pn_i(G)$ neighbours
in~$G_i$ (Lemma~\ref{lem:deg_into_cls}), (\textit{c})~the~$G_i$ are
all connected (Lemma~\ref{lem:G_i_conn}), and (\textit{d})~any
induced balanced bipartite subgraph of~$G$ with high minimum degree
contains a \emph{$k$-factor} (Lemma~\ref{lem:k-factor}),
i.e., a \(k\)-regular spanning subgraph.

\subsection{Preliminaries}
\label{sec:preliminaries}

We shall use the following Chernoff bound
(see~\cite[Corollary~A.1.14]{alon16:_probab_method}
or~\cite[Theorem~2.8]{Janson2000random}).

\begin{lemma}[Chernoff bound]
  \label{lem:chernoff}
  Let $X_1,\dots,X_n$ be independent Bernoulli random variables with
  $\PP(X_i=1) = p_i$ for all~$i$. Let $X = \sum_{i=1}^n X_i$ and
  $\mu=\mathbb{E}(X) = \sum_{i=1}^n p_i $. Then, for every
  $\varepsilon > 0$, there is $c_\varepsilon > 0$ such that
  $\PP(|X-\mu|>\varepsilon\mu) \leq 2e^{-c_\varepsilon\mu}$.
\end{lemma}

The next lemma asserts that a binomial random variable~$\Bin(n,p)$
with parameters~$n$ and~$p$ is well distributed among the congruence
classes~${}\!\bmod k$ as long as~$np(1-p)\to\infty$.

\begin{lemma}\label{lem:binomialmodk}
  Let $n\geq1$, $0 < p < 1$ and~\(k\geq2\) be given.  There is a
  positive constant $a_k$ that depends only on~$k$ 
  such that, for any integer~$t$,
  \begin{equation}
    \label{eq:binomialmodk}
    \left|\PP\big(\Bin(n,p)\equiv t\pmod* k\big)-\frac{1}{k}\right| 
    \leq e^{-a_knp(1-p)}.
  \end{equation}
\end{lemma}
\begin{proof}
  Let $\alpha$ be a primitive $k$th root of unity. 
  The $k$th roots of unity are precisely $1=\alpha^0$, $\alpha^1$, $\dots$, $\alpha^{k-1}$. 
  Now, for every $j\in\{0,\ldots,k-1\}$, we have
  \begin{equation}\label{lem:binomialmodk-eq1}
    \alpha^{-jt}\big(1+(\alpha^j-1)p\big)^n
    =\sum_{i=0}^n\binom{n}{i}p^i(1-p)^{n-i}\alpha^{j(i-t)}.
  \end{equation}
  Fix any~$t\in\{0,\dots,k-1\}$.  For any integer \(i\),
  \begin{equation*}
    \sum_{j=0}^{k-1}(\alpha^j)^{i-t}
    =
    \begin{cases}
      (\alpha^{k(i-t)}-1)/(\alpha^{i-t}-1)=0  &\text{if }k\nmid i-t,\\
      k &\text{otherwise}.
    \end{cases}
  \end{equation*}
  Thus, summing~\eqref{lem:binomialmodk-eq1} for $j=0,\dots,k-1$, we
  get
  \begin{equation}
    \label{eq:2}
    \begin{split}
      \sum_{j=0}^{k-1}\alpha^{-jt}(1+(\alpha^j-1)p)^n
      &=k\cdot\sum_{i\equiv t\pmod*k}\binom{n}{i}p^i(1-p)^{n-i}\\ 
      &=k\cdot\mathbb{P}\big(\Bin(n,p)\equiv t\pmod*k\big).
    \end{split}
  \end{equation}
  Note that, if \(j=0\), then \(\alpha^{-jt}(1+(\alpha^j-1)p)^n=1\), 
  and if $j\neq 0$, then
  \begin{equation}
    \label{eq:3}
    \begin{split}
      |(1+(\alpha^j-1)p)^n|
      &=|(1-p)+\alpha^jp|^n \\
      &=(1-2(1-\Re(\alpha^j))p(1-p))^{n/2} \\
      &\leq e^{-(1-\Re(\alpha^j))p(1-p)n}.
    \end{split}
  \end{equation}
  Inequality~\eqref{eq:binomialmodk} follows from~\eqref{eq:2}
  and~\eqref{eq:3}. 
\end{proof}

Given a graph~$G$ and~$U$, $W\subset V(G)$ with~$U\cap W=\emptyset$,
let~$e(U,W)=e_G(U,W)$ be the number of edges in~$G$ with one endpoint
in~$U$ and the other endpoint in~$W$.

\begin{definition}[$(p,\alpha)$-bijumbled]
  \label{def:(p,alpha)-bijumbled}
  Let~$p$ and~$\alpha$ be given.  We say that a graph~$G$ is
  \textit{weakly $(p,\alpha)$-bijumbled} if, for all~$U$,
  $W\subset V(G)$ with $U\cap W=\emptyset$
  and~$1\leq|U|\leq|W|\leq pn|U|$, we have
  \begin{equation}
    \label{eq:(p,alpha)-bijumbled_def}
    \big|e(U,W)-p|U||W|\big|\leq\alpha\sqrt{|U||W|}.
  \end{equation}
  If~\eqref{eq:(p,alpha)-bijumbled_def} holds for all pairs of
  disjoint sets~$U$, $W\subset V(G)$, then we say that~$G$ is
  \textit{$(p,\alpha)$-bijumbled}.
\end{definition}

\begin{fact}
  \label{fact:jumbled}
  If~$G$ is weakly $(p,\alpha)$-bijumbled, then for
  every~$U\subset V(G)$ we have
  \begin{equation}
    \label{eq:(p,alpha)-single}
    \left|
      e(G[U]) - p \binom{|U|}{2}
    \right| \leq \alpha|U|.
  \end{equation}
\end{fact}
\begin{proof}[Proof sketch]
  Let~$u=|U|$.  Double counting shows that
  \begin{equation}
    \label{eq:double_c}
    2e(G[U]){u-2\choose\lfloor u/2\rfloor-1}
    =\sum_{U'}e_G(U',U\setminus U'),
  \end{equation}
  where the sum ranges over~$U'\subset U$
  with~$|U'|=\lfloor u/2\rfloor$.
  Inequality~\eqref{eq:(p,alpha)-single} follows
  from~\eqref{eq:double_c}.  We omit the details.
\end{proof}

\begin{lemma}[Lemma~3.8 in~\cite{haxell95ramsey}]\label{lem:bijumbled}
  For any~$0<p=p(n)\leq1$, the random graph~$G(n,p)$ is a.a.s.\ weakly
  $(p,A\sqrt{pn})$-bijumbled for a certain absolute constant~$A$.  
\end{lemma}

In~\cite{haxell95ramsey}, Lemma~\ref{lem:bijumbled} is proved
with~$A=e^2\sqrt6$. 

\begin{cor}\label{cor:exists_edge}
  Let~$G=G(n,p)$, where~$0<p=p(n)\leq1$.  If $t^2>A^2n/p$, then
  a.a.s.~$e(U,W)>0$ for any pair of sets $U$, $W\subset V(G(n,p))$
  with $U\cap W=\emptyset$ and $\min\{|U|,|W|\}\geq t$.
\end{cor}

\begin{cor}\label{cor:G(n,p)-bijumbled.2}
  Suppose~$pn\geq C\log n$ for some constant~$C>3$.  Then a.a.s.~$G(n,p)$ is
  $(p,A\sqrt{pn})$-bijumbled for some~$A\geq2$. 
\end{cor}
\begin{proof}[Proof sketch]
  Lemma~\ref{lem:bijumbled} tells us that~$G(n,p)$ is a.a.s.\
  weakly $(p,A\sqrt{pn})$-bijumbled for some~$A$.  We may assume
  that~$A\geq2$.  Now let~$U$ and~$W$ be disjoint, with~$|W|>pn|U|$.
  Then~$A\sqrt{pn|U||W|}>Apn|U|$.  In particular,
  $p|U||W|-A\sqrt{pn|U||W|}\leq p|U|n-Apn|U|\leq0\leq e(U,W)$.
  As~$np\geq C\log n$ and~$C>3$, we have that~$\Delta(G(n,p))\leq2pn$
  almost surely.  Therefore
  $e(U,W)\leq2pn|U|\leq Apn|U|\leq p|U||W|+A\sqrt{pn|U||W|}$, and we
  conclude that~$G(n,p)$ is indeed $(p,A\sqrt{pn})$-bijumbled.
\end{proof}

We shall also need the following fact.  A result as in the lemma below
can be proved by considering balanced bipartitions chosen uniformly at
random and by applying a Chernoff bound for hypergeometric
distributions, but we give the result below, which can be proved by
applying Lemma~\ref{lem:chernoff} ($c_{1/5}$~below is the constant
given by Lemma~\ref{lem:chernoff} for~$\epsilon=1/5$).

\begin{lemma}
  \label{lem:split2}
  Let~$J$ be a graph on~$2u\leq n$ vertices and suppose that
  $\delta(J)\geq 10c_{1/5}^{-1}(\log n+\omega)$,
  where~$\omega=\omega(n)\to\infty$ as~$n\to\infty$.  Then, if~$n$ is
  large enough, there is $U\subset V(J)$ with~$|U|=u$ such that the
  bipartite graph $J[U,W]$ induced between~$U$ and~$W=V(J)\setminus U$
  is such that $\delta(J[U,W])\geq2\delta(J)/5$.
\end{lemma}
\begin{proof}
  Let \(\{x_1,y_1\}, \ldots, \{x_u,y_u\}\) be an arbitrary partition
  of~\(V(J)\) into pairs, and let \(U=\{z_i:1\leq i\leq u\}\), where
  each~$z_i$ is chosen uniformly at random from \(\{x_i,y_i\}\),
  independently for each~$i$.  Let \(W = V(J)\setminus U\) and put
  \(J' = J[U,W]\).  For each \(v\in V(J)\), let~\(P_v\) be the
  number of pairs~$\{x_i,y_i\}$ contained entirely in~\(N_J(v)\) and
  let~$Q_v$ be the number of~$\{x_i,y_i\}$ with
  $|\{x_i,y_i\}\cap N_J(v)|=1$.  Clearly, $\deg_J(v)=2P_v+Q_v$.  Let
  \(A = \{v : P_v < 2\deg_J(v)/5\}\).  By the definition of~\(A\), if
  \(v\notin A\), then \(\deg_{J'}(v) \geq 2\deg_J(v)/5\).  In what
  follows, we deal with the vertices in~\(A\).  Fix a vertex
  \(v\in A\).  A moment's thought tells us that
  \(\deg_{J'}(v) = P_v + d'(v)\), where~\(d'(v)\sim\Bin(Q_v,1/2)\).
  Let~$\mu=\EE(d'(v))=Q_v/2$.  Since~$v\in A$, we
  have~$\mu=(\deg_J(v)-2P_v)/2>\deg_J(v)/10\geq\delta(J)/10$.  We have
  \begin{multline*}
    \PP\[\deg_{J'}(v)<{2\over5}\deg_J(v)\]
    =\PP\[P_v+d'(v)<{2\over5}\deg_J(v)\]
    =\PP\[d'(v)<{2\over5}\deg_J(v)-P_v\]\\
    \leq\PP\[d'(v)<{2\over5}(\deg_J(v)-2P_v)\]
    =\PP\[d'(v)<{2\over5}Q_v\]
    \leq\PP\[\big|d'(v)-\mu\big|>{1\over5}\mu\],
  \end{multline*}
  which, by Lemma~\ref{lem:chernoff} and our hypothesis
  on~$\delta(J)$, is at most
  $2e^{-c_{1/5}\mu}=2e^{-c_{1/5}\delta(J)/10}=o(n^{-1})$.  Thus, by
  the union bound, the probability that $\deg_{J'}(v)<(2/5)\deg_J(v)$
  for some~$v\in A$ is~$o(1)$, showing that, for large~$n$, most
  choices of~$U$ will do (recall that the vertices $ v\notin A$ are
  never a problem).
\end{proof}

\subsection{Degree classes of~$G(n,p)$}
\label{sec:degree-classes-gnp}

We first show that the degree classes of~$G(n,p)$ are typically of
cardinality about~$n/k$.  This is assertion~(\textit{a}) given at the
beginning of Section~\ref{sec:technical-lemmas}.  In
Lemmas~\ref{lem:deg_into_cls}, \ref{lem:G_i_conn}
and~\ref{lem:k-factor}, we prove assertions~(\textit{b}), (\textit{c})
and~(\textit{d}).

\begin{lemma}\label{lem:bal_deg_cls}
  Let~$k\geq2$ be a fixed integer and let~$p=p(n)$
  with~$np(1-p)\to\infty$ as~$n\to\infty$ be given.  Then, with
  probability at least~$1-o(1/n)$, for every~$1\leq i\leq k$ we have
  \begin{equation}
    \label{eq:bal_deg_cls}
    {n\over2k}\leq n_i(G(n,p))\leq{3n\over2k}.
  \end{equation}
\end{lemma}
\begin{proof}
  Fix~$i$ ($1\leq i\leq k$).  We show that~\eqref{eq:bal_deg_cls}
  holds with probability~$1-o(1/n)$.  The result then follows from the
  union bound.
  
  Let~$G=G(n,p)$.  Fix~$U\subset V(G)$
  with~$|U|=\lceil(1-1/4k)n\rceil$ and let~$W=V(G)\setminus U$.
  Let~$m=|W|$.  Let~$F\subset{U\choose2}$ and condition on~$E(H)=F$,
  where~$H=G[U]$.  For every~$u\in U$, let~$X_u$ be the indicator
  function of the event $\{\deg_G(u)\equiv i\pmod*k\}$.  Since
  $\deg_G(u)=\deg_H(u)+e_G(\{u\},W)$ and we are conditioning
  on~$E(H)=F$, we have that
  $p_u=\PP(X_u=1)=\PP\big(\Bin(m,p)\equiv t\pmod*k\big)$,
  where~$t=i-\deg_H(u)$.  Lemma~\ref{lem:binomialmodk} tells us that
  $|p_u-1/k|\leq e^{-a_kp(1-p)m}=o(1)$.  Let~$X=\sum_{u\in U}X_u$ and
  note that~$\EE(X)=|U|(1/k+o(1))$.  Lemma~\ref{lem:chernoff} then
  tells us that, for some absolute constant~$c>0$, 
  \begin{equation*}
    \PP\left(\left|X-{1\over k}|U|\right|>{n\over4k}\right)
    \leq2e^{-c|U|/k}=o\left(1\over n\right).
  \end{equation*}
  Also, note that $X\leq n_i(G)\leq X+m$, and that
  \begin{equation*}
    {1\over k}|U|-{n\over4k}\geq{1\over
      k}\llceil\[1-{1\over4k}\]n\rrceil-{n\over 4k}
    \geq{n\over2k}
  \end{equation*}
  and
  \begin{equation*}
    {1\over k}|U|+{n\over4k}+m\leq{1\over
      k}\llceil\[1-{1\over4k}\]n\rrceil+{n\over 4k}+{n\over4k}
    \leq{3n\over2k}.
  \end{equation*}
  Therefore
  \begin{equation*}
    \PP\[{n\over2k}\leq n_i(G)\leq{3n\over2k}\;\bigg|\;E(H)=F\]
    =1-o\[1\over n\]. 
  \end{equation*}
  Since this holds for arbitrary~$F$, the result follows.
\end{proof}

Recall that \(V_i=V_i(G)\) is the set of vertices~$v$ of~\(G\) with
$\deg_G(v)\equiv i\pmod*k$ and $G_i=G[V_i]$.  

\begin{lemma}\label{lem:deg_into_cls}
  For every integer~$k\geq2$ there is a positive constant~$C$ such
  that if~$p=p(n)\geq Cn^{-1}\log n$ and $n(1-p)\to\infty$ as
  $n\to\infty$, then a.a.s.~$G=G(n,p)$ is such that,
  for every~$v\in V(G)$ and every~$1\leq i\leq k$,
  \begin{equation*}
    |N(v)\cap V_i|\geq{pn\over3k}.
  \end{equation*}
\end{lemma}
\begin{proof}
  Let~$c_{1/4}>0$ be as given by Lemma~\ref{lem:chernoff} and
  let~$C=3k/c_{1/4}$.  We prove that this choice of~$C$ will do.
  Fix~$1\leq i\leq k$ and~$v\in V=V(G)$.  Let~$U=V\setminus\{v\}$.  We
  first generate the edges of~$G=G(n,p)$ in~$H=G[U]$.
  Since~$p\geq Cn^{-1}\log n$, our assumption that~$n(1-p)\to\infty$
  implies that~$np(1-p)\to\infty$ as well.  Hence
  Lemma~\ref{lem:bal_deg_cls} applies and we see that, with
  probability $1-o(1/n)$, we have
  \begin{equation}
    \label{eq:5}
    n_j(H)\geq{n-1\over2k}
  \end{equation}
  for all $1\leq j\leq k$.  Let us suppose that~\eqref{eq:5} does hold
  for every~$j$.  We now generate the edges between~$v$ and~$U$ in~$G$.
  Clearly, $N(v)\cap V_i=N(v)\cap V_{i-1}(H)$, where, of course, we
  consider the indices modulo~$k$.  Also,
  $|N(v)\cap V_i|\sim\Bin(n_{i-1}(H),p)$.  Note that
  \begin{equation*}
    \EE(|N(v)\cap V_i|)=pn_{i-1}(H)\geq p\,{n-1\over2k}
    \geq{4pn\over9k}
  \end{equation*}
  and also that $(3/4)pn_{i-1}(H)\geq{pn/3k}$ for all large
  enough~$n$.  Lemma~\ref{lem:chernoff} then gives that, with
  probability $1-2\exp(-4c_{1/4}pn/9k)=1-o(1/n)$, we have
  $|N(v)\cap V_i|\geq(3/4)pn_{i-1}(H)\geq pn/3k$.  It now suffices
  to take the union bound considering all~$1\leq i\leq k$
  and~$v\in V$.
\end{proof}

\begin{lemma}\label{lem:G_i_conn}
  Let~$k\geq2$ be an integer and let~$C$ and~$p=p(n)$ be as in
  Lemma~\ref{lem:deg_into_cls}.  Then a.a.s.~$G=G(n,p)$
  is such that~$G_i$ is connected for every~$1\leq i\leq k$.
\end{lemma}
\begin{proof}
  Fix~$i$.  Lemma~\ref{lem:deg_into_cls} tells us that a.a.s.
  \begin{equation}
    \label{eq:6}
    \delta(G_i)\geq{pn\over3k}.
  \end{equation}
  We suppose~\eqref{eq:6} holds and that~$G$ is
  $(p,A\sqrt{pn})$-bijumbled for some~$A\geq2$ (recall
  Corollary~\ref{cor:G(n,p)-bijumbled.2}) and deduce that~$G_i$ is
  connected if~$n$ is large enough.  Suppose for a contradiction
  that~$J$ is a component of~$G_i$ with~$t=|V(J)|\leq|V(G_i)|/2$.  The
  number~$e(J)$ of edges in~$J$ satisfies
  \begin{equation*}
    {1\over2}t\,{pn\over3k}\leq{1\over2}t\delta(J)\leq e(J)
     \stackrel{\text{Fact~\ref{fact:jumbled}}}{\leq}
     p{t\choose2}+A\sqrt{pn}\,t
    \leq p{t^2\over2}+A\sqrt{pn}\,t, 
  \end{equation*}
  whence
  \begin{equation}
    \label{eq:8}
    {pn\over6k}\leq{1\over2}pt+A\sqrt{pn}.
  \end{equation}
  Since~$pn\to\infty$, it follows from~\eqref{eq:8} that, say,
  $t\geq2n/7k$ for any large enough~$n$.  By the choice of~$J$, we
  have $|V(G_i)\setminus V(J)|\geq t$.  Therefore, by
  Corollary~\ref{cor:exists_edge}, we have
  $e(V(J),V(G_i)\setminus V(J))>0$, as
  $t^2\gg A^2n^2/C\log n\geq A^2n/p$.  Since~$J$ is a component
  of~$G_i=G[V_i]$ this is a contradiction.  We conclude that~$G_i$ is
  indeed connected.
\end{proof}

\begin{lemma}\label{lem:k-factor}
  Let~$k\geq1$, $c>0$, $A>0$ and~$0<p=p(n)<1$ be given.  Suppose~$G$
  is a $(p,A\sqrt{pn})$-bijumbled graph of order~$n$ and~$p\gg1/n$.
  Then, if~$n$ is large enough, for any~$U$ and~$W\subset V(G)$ with
  $U\cap W=\emptyset$, $|U|=|W|\geq cn$ and
  $\delta(G[U,W])\geq p|U|/8$, the graph~$G[U,W]$ contains a
  $k$-factor.
\end{lemma}
\begin{proof}
  Let~$U$ and~$W$ be as in the statement of the lemma and
  let~$m=|U|=|W|$.  We prove that~$G[U,W]$ contains a $k$-factor by
  induction on~$k$.  Fix~$k\geq1$ and suppose~$G[U,W]$ contains a
  $(k-1)$-factor~$F$.  It suffices to prove that~$B=G[U,W]-F$ contains
  a perfect matching.  We check Hall's condition: for
  every~$S\subset U$, we have $|N(S)|\geq|S|$.
  Let~$\delta=\delta(B)$.  A simple argument shows that if
  $|S|\leq\delta$ or $|S|>m-\delta$, then $|N(S)|\geq|S|$.  We
  therefore assume that $\delta<|S|\leq m-\delta$ and suppose for a
  contradiction that $|N(S)|<|S|$.  Let~$s=|S|$.  Then
  \begin{equation*}
    \[{1\over8}cpn-k\]s\leq\[{1\over8}pm-k\]s\leq\delta s
    \leq e(S,N(S))\leq ps^2+A\sqrt{pn}\,s,
  \end{equation*}
  whence
  \begin{equation}
    \label{eq:7}
    {1\over8}cpn-k\leq ps+A\sqrt{pn}.
  \end{equation}
  Since~$pn\gg1$, it follows from~\eqref{eq:7} that if~$n$ is large
  enough, then, say, $cpn/9\leq ps$ and hence~$|S|=s\geq cn/9$.

  Let~$T=W\setminus N(S)$ and note that $N(T)\subset U\setminus S$.
  Hence, $|N(T)|\leq|U|-|S|<|U|-|N(S)|=|W|-|N(S)|=|T|$.  Arguing as
  above, we get that $|T|\geq cn/9$.  Using that $pn\gg1$, we see by
  Corollary~\ref{cor:exists_edge} that $e(S,T)>0$, which contradicts
  the definition of~$T$.  We conclude that $B=G[U,W]-F$ satisfies
  Hall's condition.  This concludes the induction step and the result
  follows.
\end{proof}

\section{Main theorem}
\label{sec:main-theorem}
	
We now prove Theorem~\ref{thm:main}.  In what follows, we say that a
graph~\(G\) is a \emph{mod~$k$ graph} if all its non-isolated vertices
have degrees congruent to~$1$ mod~$k$.  The reader may find it useful
to recall the notation and terminology introduced at the beginning of
Section~\ref{sec:technical-lemmas}.
		
The key idea of this proof when $k$ and $n$ are even or $k$ is odd is
(1)~to use Lemma~\ref{lem:deg_into_cls} to find a set of
vertex-disjoint stars forming a star forest~$F$
with~$\chi_k'(F)\leq k-1$ so that~$G'=G-E(F)$ has an even number of
vertices in each degree class~\(V_i(G')\) with \(i>1\), and then
(2)~use Lemma~\ref{lem:k-factor} to find, for each \(i>1\), a
bipartite \((i-1)\)-factor~$B_i$ in~\(G'_i\) and then let
$G''=G'-E(B)$, where $B=\bigcup_{2\leq i\leq k}B_i$ (see
Figure~\ref{fig}).  By construction, $\chi_k'(F\cup B)\leq k-1$ and
$G''=G-E(F\cup B)$ is a mod~$k$ graph, which can be colored
monochromatically.  It follows that $\chi_k'(G)\leq k$.  The remaining
case, namely when~$k$ is even and~$n$ is odd, then follows by using
the $n$~even case to color $G-v$ with~$k$ colors for some~$v\in V(G)$,
and then coloring most of the edges incident to \(v\) with
the~\((k+1)\)st color.
	
\begin{proof}[Proof of Theorem~\ref{thm:main}]
  Let~$C_k=\max\{C,41c_{1/5}^{-1}k\}$, where~$C$ is the constant given by
  Lemma~\ref{lem:deg_into_cls} and $c_{1/5}$~is the constant given by
Lemma~\ref{lem:chernoff} for~$\epsilon=1/5$.  Let~$p=p(n)$ be as in the statement
  of the theorem and let~$G=G(n,p)$.  Below, we tacitly assume
  that~$n$ is large enough whenever necessary.

  We start by observing that $\chi'_k(G)\geq k$ holds a.a.s.\
  regardless of the parity of~$k$.  
  Indeed, owing to our hypothesis on~$p$, Lemma~\ref{lem:bal_deg_cls}
  tells us that \(n_k(G(n,p))\geq n/2k\) with probability at least \(1 - o(1/n)\).
  Noting that~$G$ a.a.s.\ has no isolated vertices, we deduce that~$G$
  a.a.s.\ has a vertex~$v$ of nonzero degree with
  $\deg_G(v)\equiv0\pmod*k$.  It is clear that such a vertex~$v$
  forces $\chi'_k(G)\geq k$, regardless of the parity of~$k$.  We also
  have to prove that, for even~$k$ and odd~$n$, we a.a.s.\ have
  $\chi'_k(G)\geq k+1$.  This is done below.

  We divide the remainder of the proof according to the parity
  of~\(k\).

  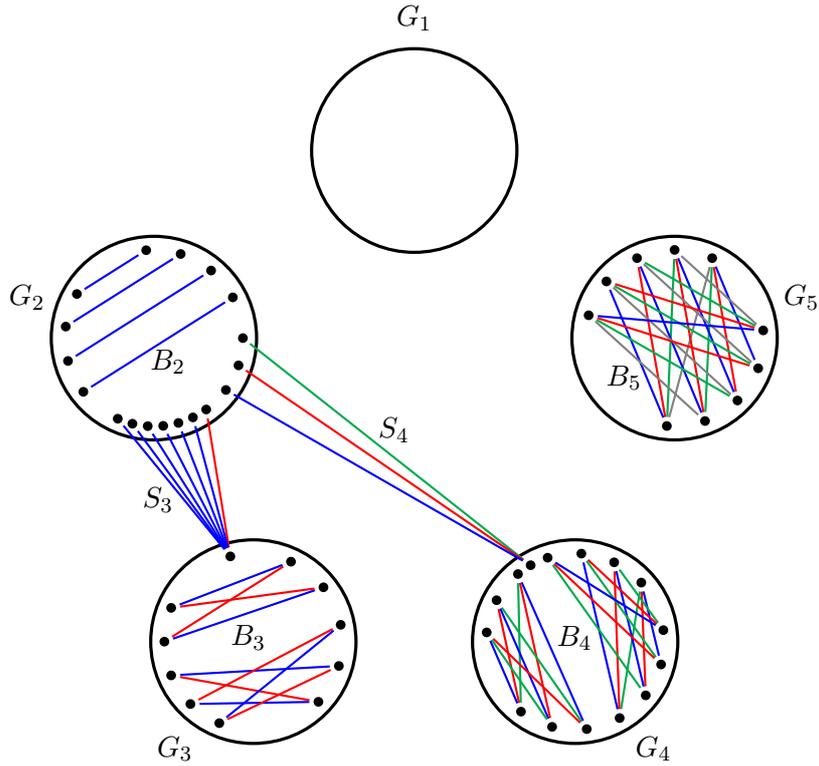
\begin{figure}[!ht]
    \centering
    \begin{tikzpicture}[scale=.9, thick]

  \tikzset{black vertex/.style={circle,draw,minimum size=1mm,inner sep=0pt,outer sep=2pt,fill=black, color=black}}

  \tikzstyle{every circle node} = [draw]

  \foreach \i in {1,...,5}
  {
    \node (\i) [] at (\i*72+18:4) {};
    \draw[very thick] (\i) circle (1.5cm);
    \node (label) [] at (\i*72+18:5.95) {$G_\i$};
  }
  
  \node (21) [black vertex] at ($(2)+(0:1.3)$) {};
  \node (22) [black vertex] at ($(2)+(-18:1.3)$) {};
  \node (23) [black vertex] at ($(2)+(-36:1.3)$) {};
  \node (24) [black vertex] at ($(2)+(-54:1.3)$) {};
  \node (25) [black vertex] at ($(2)+(-64:1.3)$) {};
  \node (26) [black vertex] at ($(2)+(-74:1.3)$) {};
  \node (27) [black vertex] at ($(2)+(-84:1.3)$) {};
  \node (28) [black vertex] at ($(2)+(-94:1.3)$) {};
  \node (29) [black vertex] at ($(2)+(-104:1.3)$) {};
  \node (210) [black vertex] at ($(2)+(-114:1.3)$) {};

  \node (31) [black vertex] at ($(3)+(105:1.3)$) {};

  \draw[red, thick] (31) -- (24);
  \begin{scope}[blue]
    \draw (31) -- (25);
    \draw (31) -- (26);
    \draw (31) -- (27);
    \draw (31) -- (28);
    \draw (31) -- (29);
    \draw (31) -- (210) node[black, pos=0.4, left]{$S_3$};
  \end{scope}
  
  \node (41) [black vertex] at ($(4)+(120:1.3)$) {};

  \draw[yellow!10!green] (41) -- (21) node[black, pos=0.6, right]{$\;S_4$};
  \draw[red] (41) -- (22);
  \draw[blue] (41) -- (23);      

  \foreach \i in {1,...,4} {
    \node (2a\i) [black vertex] at ($(2)+(5+22.5*\i:1.3)$) {};
    \node (2b\i) [black vertex] at ($(2)+(240-22.5*\i:1.3)$) {};
  }
  \draw[blue] (2a1) -- (2b1) node[black, pos=0.7, right]{$\;B_2$};
  \draw[blue] (2a2) -- (2b2);
  \draw[blue] (2a3) -- (2b3);
  \draw[blue] (2a4) -- (2b4);

  \foreach \i in {1,...,5}{
    \node (3a\i) [black vertex] at ($(3)+(30-100+27*\i:1.3)$) {};
    \node (3b\i) [black vertex] at ($(3)+(30+240-22.5*\i:1.3)$) {};
  }
  \begin{scope}[thick]
    \draw[blue] (3a1) -- (3b2);
    \draw[red] (3a1) -- (3b3);
    \draw[red] (3a2) -- (3b1);
    \draw[blue] (3a2) -- (3b3);
    \draw[blue] (3a3) -- (3b1);
    \draw[red] (3a3) -- (3b2);
    \draw[blue] (3a4) -- (3b4)  node[black, pos=0.5, below]{$\;B_3$};
    \draw[red] (3a4) -- (3b5);
    \draw[red] (3a5) -- (3b4);
    \draw[blue] (3a5) -- (3b5);
  \end{scope}

  \foreach \i in {1,2,3,4,5,6,7}{
    \node (4a\i) [black vertex] at ($(4)+(20+22*\i:1.3)$) {};
    \node (4b\i) [black vertex] at ($(4)+(30-22.5*\i:1.3)$) {};
  }
  \draw[blue] (4a1) -- (4b2);
  \draw[red] (4a1) -- (4b3);
  \draw[yellow!10!green] (4a1) -- (4b4);
  \draw[yellow!10!green] (4a2) -- (4b1);
  \draw[blue] (4a2) -- (4b3);
  \draw[red] (4a2) -- (4b4);
  \draw[red] (4a3) -- (4b1);
  \draw[yellow!10!green] (4a3) -- (4b2);
  \draw[blue] (4a3) -- (4b4);
  \draw[blue] (4a4) -- (4b1);
  \draw[red] (4a4) -- (4b2);
  \draw[yellow!10!green] (4a4) -- (4b3);

  \draw[blue] (4a5) -- (4b5) node[black, pos=0.4, right]{$B_4$};
  \draw[red] (4a5) -- (4b6);
  \draw[yellow!10!green] (4a5) -- (4b7);
  \draw[yellow!10!green] (4a6) -- (4b5);
  \draw[blue] (4a6) -- (4b6);
  \draw[red] (4a6) -- (4b7);
  \draw[red] (4a7) -- (4b5);
  \draw[yellow!10!green] (4a7) -- (4b6);
  \draw[blue] (4a7) -- (4b7);

  \foreach \i in {1,2,3,4,5}{
    \node (5a\i) [black vertex] at ($(5)+(40+25*\i:1.3)$) {};
    \node (5b\i) [black vertex] at ($(5)+(30-25*\i:1.3)$) {};
  }
  \draw[blue] (5a1) -- (5b2);
  \draw[red] (5a1) -- (5b3);
  \draw[yellow!10!green] (5a1) -- (5b4);
  \draw[gray] (5a1) -- (5b5);
  \draw[gray] (5a2) -- (5b1);
  \draw[blue] (5a2) -- (5b3);
  \draw[red] (5a2) -- (5b4);
  \draw[yellow!10!green] (5a2) -- (5b5);
  \draw[yellow!10!green] (5a3) -- (5b1);
  \draw[gray] (5a3) -- (5b2);
  \draw[blue] (5a3) -- (5b4);
  \draw[red] (5a3) -- (5b5);
  \draw[red] (5a4) -- (5b1);
  \draw[yellow!10!green] (5a4) -- (5b2);
  \draw[gray] (5a4) -- (5b3);
  \draw[blue] (5a4) -- (5b5);
  \draw[blue] (5a5) -- (5b1);
  \draw[red] (5a5) -- (5b2);
  \draw[yellow!10!green] (5a5) -- (5b3);
  \draw[gray] (5a5) -- (5b4) node[black, pos=0.6, left]{$B_5\,$};
  \end{tikzpicture}
    \caption{%
      The subgraphs~$F=S_3\cup S_4$ and~$B=B_2\cup\cdots\cup B_5$.
      Here, $k=5$, $n_2(G)$ and~$n_5(G)$ are even and $n_3(G)$
      and~$n_4(G)$ are odd.}
    \label{fig}
  \end{figure}
		
  \noindent
  \textbf{Case 1. $k$ even.} We consider the even~$n$ and odd~$n$
  cases separately.
  
  \smallskip
  \noindent
  \textbf{Case 1.1 $n$ even.}  We shall prove that, in this case,
  $\chi_k'(G)\leq k$ holds a.a.s.\ by following the strategy outlined
  at the beginning of Section~\ref{sec:main-theorem}.

  We first claim that, a.a.s., there are vertex-disjoint stars
  \(S_3,\dots,S_k\) such that, for each $i\in\{3,4,\dots,k\}$, the
  star~\(S_i\) has \(i-1\) edges, is centered in~$G_i$ and its leaves
  belong to~$G_2$.  This can be seen by applying
  Lemma~\ref{lem:deg_into_cls} successively, to obtain each of
  the~$S_i$ ($i=3,4,\dots$) in turn.  Let $F=\bigcup_{i\in I}S_i$,
  where $I=\{i\:n_i(G)=|V(G_i)|\text{ is odd}\}$.

  Let~$G'=G-E(F)$.  Note that the vertices of~$G$ incident to the
  edges of~$F$ have their degrees changed by the removal of~$F$: they
  all become of degree $1\pmod*k$ in~$G'$, and therefore each of them
  ``moves'' from~$G_i$ for some~$i\geq2$ to~$G_1'$.  Furthermore, note
  that~$n_i(G')=|V(G_i')|$ is even for all~$i\geq 3$,
  because~$G_i'=G_i$ if~$n_i(G)$ is even and~$G_i'=G_i-v_i$
  if~$n_i(G)$ is odd, where~$v_i$ is the center of the star~$S_i$.
  
  Since~$n$ is even, the number~$n_{\rm even}'$ of vertices of even
  degree in~$G'$ has the same parity as~$n_{\rm odd}'$, the number of
  vertices of odd degree in~$G'$, which is even.  As~$k$ is even,
  $n_{\rm even}'=\sum\{n_i(G')\:i\text{ even}\}$ and the fact
  that~$n_i(G')$ is even for all~$i\geq3$ implies that~$n_{\rm even}'$
  and~$n_2(G')$ have the same parity, and hence~$n_2(G')$ is even.  We
  conclude that $n_i(G')$ is even for all~$i\geq2$.

  Fix~$2\leq i\leq k$.  In view of Lemmas~\ref{lem:bal_deg_cls}
  and~\ref{lem:deg_into_cls}, we assume that~\eqref{eq:bal_deg_cls}
  holds and that we have $\delta(G_i)\geq pn/3k$, and
  consequently~$\delta(G_i')\geq\delta(G_i)-1\geq
  pn/4k\geq(C_k/4k)\log n\geq(41/4c_{1/5})\log n$.  We now apply
  Lemma~\ref{lem:split2} with~$J=G_i'$ and obtain~$U$
  and~$W\subset V(G_i')$
  with~$|U|=|W|=n_i(G')/2\leq n_i(G)/2\leq3n/4k$
  and~$U\cap W=\emptyset$ such that
  \begin{equation*}
    \delta(G[U,W])=\delta(G_i'[U,W])
    \geq{2\over5}\delta(G_i')\geq{pn\over10k}\geq{1\over8}p|U|.
  \end{equation*}
  Note that
  $|U|=|W|=n_i(G')/2\geq(n(G_i)-1)/2\geq(n/2k-1)/2\geq n/5k$.  We are
  now in a position to apply Lemma~\ref{lem:k-factor} and obtain an
  $(i-1)$-factor~$B_i$ in~$G[U,W]=G_i'[U,W]$.

  Let~$B=\bigcup_{2\leq i\leq k}B_i$ and~$G''=G'-E(B)$.  Note
  that~$G''$ is a mod~$k$ graph, which can be entirely colored with
  color~$1$, say.  Furthermore, $\chi_k'(F\cup B)\leq k-1$.  We
  conclude that~$\chi_k'(G)\leq k$.  This finishes the proof of
  Case~1.1.
  
  Before we proceed, for later reference, we observe the following:
  since \(G''\) contains every edge of \(G\) incident to the vertices
  in \(G_1\), the coloring we have obtained in this case is such that
  every vertex in~$G_1$ is incident only to edges of a certain fixed
  color (color~\(1\) above).
		
  \smallskip
  \noindent
  \textbf{Case 1.2 $n$ odd.} Fix a vertex $u$ in~$G=G(n,p)$ and
  let~$H=G-u$.  Lemma~\ref{lem:deg_into_cls} tells us that we may
  suppose that~$u$ has at least~$pn/3k$ neighbors in~$H_1$.

  By Case 1.1, with high probability~$H$ can be colored with colors
  $1,\dots,k$ so that all edges incident to vertices with degree
  \(1 \pmod* k\) in \(H\) are colored with the same color, say \(1\)
  (see the last paragraph of Case~1.1).  We now color the edges
  incident to \(u\).  Suppose $\deg(u)\equiv d\pmod k$ where
  $d\in\{1,2,\dots,k\}$.  If \(d\neq 1\), then we assign each of the
  colors $2,\dots,d$ once to an edge joining $u$ to vertices in $H_1$
  (this is possible since there are at least $k-1$ such edges),
  leaving a number congruent to~\(1\pmod*k\) of uncolored edges
  incident to \(u\).  We assign to these uncolored edges a new color.
  We thus obtain a $\chi_k'$-coloring of~$G$ with $k+1$ colors.
  
  Suppose now that \(G\) admits a $\chi_k'$-coloring with~$k$
  colors.  This implies that all edges incident to any given vertex in
  $G_1$ must get the same color.  By Lemma~\ref{lem:G_i_conn}, the
  graph~\(G_1\) is connected with high probability, and hence
  a.a.s.~all the edges of~\(G\) incident to vertices of~\(G_1\) must
  be colored with the same color, say~\(1\).  Moreover, by
  Lemma~\ref{lem:deg_into_cls}, the set~$V(G_1)$ is a.a.s.\ a
  dominating set, that is, every vertex of~$G$ not in~$V(G_1)$ is
  adjacent to some vertex in~$V(G_1)$.  This implies that a.a.s.\ the
  edges of color~\(1\) induce a spanning subgraph of~\(G\).  Let~$J$
  be this spanning subgraph.  Since~$J$ is a mod~$k$ graph and~$k$ is
  even, every vertex of~$J$ has odd degree.  This is a contradiction
  as~$J$ has~$n$ vertices and~$n$ is odd.  This argument shows that
  $\chi_k'(G)\geq k+1$ with high probability.
    
  \smallskip
  \noindent
  \textbf{Case 2. $k$ odd.}  We proceed as in Case~1.1, except that,
  to produce $G'=G-E(F)$ so that~$n_i(G')$ is even for every~$i\geq2$,
  we have to argue a little more.

  Recall~$I=\{i\:n_i(G)\text{ is odd}\}$.  If~$I\neq\emptyset$, then
  we can use the stars~$S_i$ ($i\in I$) as in Case~1.1 to
  define~$F=\bigcup_{i\in I}S_i$, except that, if doing so we
  obtain~$G'=G-E(F)$ with~$n_2(G')$ odd, then we replace the
  star~$S_i$ with $i-1$ rays by a star~$S_i'$ with $k+i-1$ rays for an
  arbitrary~$i\in I$.  Since $i-1$ and $k+i-1$ have opposite parities,
  we can thus force~$n_2(G')$ to be even.  If~$I=\emptyset$
  and~$n_2(G)$ is odd, we can take two stars~$S_3$ and~$S_3'$ with
  centers in~$V_3(G)$ and with $2$~rays and $k+2$~rays, respectively,
  and define $F=S_3\cup S_3'$.  Then~$G'=G-E(F)$ is such
  that~$n_i(G')$ is even for every~$i\geq2$.  If~$I=\emptyset$
  and~$n_2(G)$ is even, we simply let~$G'=G$.

  The rest of the proof follows Case~1.1 \textit{mutatis mutandis}.
\end{proof}
	
\section{A lower bound for $\max_G\chi_k'(G)$}
\label{sec:deterministic-lower-bound}

In this section, we present a lower bound for the maximum mod~\(k\)
chromatic index of graphs.  We clearly have $\max_G\chi_k'(G)\geq k$,
because any graph~$G$ that contains a vertex~$v$ with $\deg_G(v)>0$
and $\deg_G(v)\equiv0\pmod*k$ is such that $\chi_k'(G)\geq k$.  In
1991, Pyber~\cite{pyber1991covering} showed that his upper bound
of~\(4\) for the mod~\(2\) chromatic index of graphs is tight because
the \(4\)-wheel (the graph obtained from a cycle of length~\(4\) by
adding a new vertex adjacent to all of its vertices) has mod~\(2\)
chromatic index equal to~\(4\).  Note that the \(4\)-wheel is
precisely the complete \(3\)-partite graph~$K_{1,k,k}$ with~$k=2$.
The proposition below generalizes this observation:
\(\chi'_k(K_{1,k,k})=k+2\) for every~$k\geq2$; in particular,
$\max_G\chi_k'(G)\geq k+2$.

\begin{prop}
  \label{prop:lwbd}
  For every~$k\geq2$, we have $\chi_k'(K_{1,k,k})=k+2$.
\end{prop}
\begin{proof}
  Let \(G=K_{1,k,k}\) be the complete $3$-partite graph with vertex
  classes~$\{u\}$, $A$ and~$B$.  Suppose for a contradiction
  that~\(G\) has a $\chi_k'$-coloring with~$c$ colors, where
  $c\leq k+1$.  Note that some color, say \(1\), must be used to color
  precisely \(k+1\) edges incident to~\(u\), and hence, every other
  edge incident to~\(u\) must be colored with a distinct color.  In
  particular, this implies that \(c\geq k\).  On the other hand, given
  a vertex \(v\neq u\), there are only two ways of coloring the $k+1$
  edges incident to~\(v\): (\textit{a})~by coloring all the $k+1$
  edges with the color used on \(uv\), or (\textit{b})~by coloring
  each of the $k+1$ edges with a distinct color.  Vertices of
  type~(\textit{a}) are called \emph{monochromatic} and vertices of
  type~(\textit{b}) are called \emph{rainbow}.  Clearly, since we only
  have $k+1$ colors, every color occurs at every rainbow vertex.

  \begin{claim}
    Every vertex $v\in A\cup B$ is rainbow. 
  \end{claim}
  \begin{proof}
    Let us first note that, since there are $k+1$ edges incident
    to~$u$ with color~$1$, we may assume without loss of generality
    that there are two vertices~$x$ and~$y$ in~$A$ and a vertex~$z$
    in~$B$ for which $ux$, $uy$ and $uz$ have color~$1$.  We now
    fix~$v\in A\cup B$ and show that it is rainbow.

    \smallskip
    \noindent
    \textbf{Case 1. $uv$ has color~$1$ and $v\in B$.}  Suppose~$v$ is
    monochromatic.  Then~$v$ is monochromatic of color~$1$, as~$uv$ is
    of color~$1$.  Note that both~$x$ and~$y$ are then incident to at
    least two edges of color~$1$, and hence they are both
    monochromatic of color~$1$.  It follows that every vertex in~$B$
    is monochromatic of color~$1$. Since~$k\geq2$, this implies that
    every vertex in~$A$ is also monochromatic of color~$1$.  We
    conclude that~$u$ is also monochromatic of color~$1$, and this is
    a contradiction.  Hence~$v$ is rainbow.

    \smallskip
    \noindent
    \textbf{Case 2. $uv$ has color~$1$ and $v\in A$.}  Suppose~$v$ is
    monochromatic.  Then~$v$ is monochromatic of color~$1$.  Note
    that~$z$ is then incident to two edges of color~$1$ and hence is
    monochromatic of color~$1$.  We are now as at the beginning of
    Case~1 above (we have a vertex in~$B$ monochromatic of color~$1$),
    and hence we again have a contradiction.  Thus~$v$ must be
    rainbow.

    \smallskip
    \noindent
    \textbf{Case 3. $uv$ has a color different from~$1$ and~$v\in B$.}
    Let~$uv$ have color~$2$.  Suppose~$v$ is monochromatic.  The
    argument in Case~1 shows that~$z$ is rainbow.  The edge~$zu$ has
    color~$1$, whence there is some~$w\in A$ such that~$zw$ has
    color~2.  Since we are supposing that~$v$ is monochromatic of
    color~$2$, the edge~$vw$ is of color~$2$.  This implies that~$w$
    is monochromatic of color~$2$, giving another edge of color~$2$
    incident to~$u$.  This is a contradiction, showing that~$v$ is
    rainbow.

    \smallskip
    \noindent
    \textbf{Case 4. $uv$ has a color different from~$1$ and~$v\in A$.}
    It suffices to repeat the argument in Case~3, replacing the
    vertex~$z$ in that argument by the vertex~$x$ or~$y$.  This
    concludes the proof of the claim.
  \end{proof}
  
  Let~$2$ be a color different from~$1$ that occurs at~$u$.  We know
  that~$2$ occurs exactly once at~$u$.  Since every vertex
  in~$A\cup B$ is rainbow, color~$2$ occurs at every vertex
  in~$A\cup B$ and it clearly occurs exactly once at every such
  vertex.  This means that the edges of color~$2$ form a perfect
  matching, but this is impossible as $G=K_{1,k,k}$ has an odd number
  of vertices.  This shows that~$\chi_k'(K_{1,k,k})>k+1$.

  We now show that~$\chi_k'(K_{1,k,k})\leq k+2$.  Suppose
  $A=\{a_i\:1\leq i\leq k\}$ and $B=\{b_i\:1\leq i\leq k\}$.  Let
  $A'=A\cup\{a_{k+1}\}$ and $B'=B\cup\{b_{k+1}\}$, where~$a_{k+1}$
  and~$b_{k+1}$ are two new vertices, and consider the complete
  bipartite graph~$G^+$ with vertex classes~$A'$ and~$B'$.  Let us
  color the edges of~$G^+$ properly with colors $1,\dots,k+1$ (the
  chromatic index of~$G^+=K_{k+1,k+1}$ is $k+1$).  We now omit the
  vertices~$a_{k+1}$ and~$b_{k+1}$ from~$G^+$ and add a new vertex~$u$
  adjacent to all the vertices in~$A\cup B$.  We thus obtain
  a~$K_{1,k,k}$.  It remains to color the edges~$uv$ ($v\in A\cup B$).
  Let~$m_j$ be the `missing color' at~$b_j$ ($1\leq j\leq k$): this is
  the color of~$a_{k+1}b_j$ in the proper coloring of~$G^+$.  Note
  that all the~$m_j$ ($1\leq j\leq k$) are distinct.  We now color the
  $k+1$ edges $ua_i$ ($1\leq i\leq k$) and the edge~$ub_1$ with color
  $k+2$, and color the edges~$ub_j$ ($2\leq j\leq k$) with
  color~$m_j$.  It is then clear that the edges of color~$c$
  ($1\leq c\leq k+1$) form a matching and the edges of color~$k+2$
  form a star with $k+1$ rays.  Thus~$\chi_k'(K_{1,k,k})\leq k+2$.
\end{proof}
	
We put forward the following rather optimistic conjecture
(see~\cite[Conjeture~6]{BoCoKo20}).

\begin{conj}
  There is an absolute constant~$C$ such that $\chi'_k(G)\leq k+C$ for
  every graph $G$.
\end{conj}
	
\section{Concluding remarks and future work}
\label{sec:concluding-remarks}
	
In this paper we determined the mod~\(k\) chromatic index
of~\(G(n,p)\) for~$p=p(n)$ such that~\(p \geq C_kn^{-1}\log n\) and
\(n(1-p) \rightarrow\infty\).  It is natural to investigate the
remaining ranges of~\(p\).  For instance, if $G$ is a forest, it is
not hard to prove that
\(\chi'_k(G) = \max\big\{r\in\{1,\ldots,k\}\: r\equiv d(v)\pmod*k
\text{ for some } v \in V(G)\big\}\). This observation settles the
case $p \ll n^{-1}$, since in this range $G(n,p)$ is a.a.s.\ a
forest. The next step would be to consider $p=cn^{-1}$ for some
constant $c \in (0,1)$, in which case the components of~$G(n,p)$ are
a.a.s.\ trees and unicyclic graphs. Unfortunately, the formula above
for~$\chi'_k(G)$ does not extend to all unicyclic graphs: it is not
hard to prove that if \(G\) is any graph that contains a cycle of
length \(\ell\geq 3\) in which \(\ell-1\) vertices have degree
precisely \(k+1\), and one vertex has degree at most \(k\), then
\(\chi'_k(G)\geq k+1\).  Quite possibly, the most challenging range
would be~$n^{-1}\leq p\leq cn^{-1}\log n$, where~$c$ is a smallish
constant.
	
\bibliography{ref}

\begin{bibdiv}
\begin{biblist}

\bib{alon16:_probab_method}{book}{
      author={Alon, Noga},
      author={Spencer, Joel~H.},
       title={The probabilistic method},
     edition={4th ed.},
      series={Wiley Series in Discrete Mathematics and Optimization},
   publisher={John Wiley \& Sons, Inc., Hoboken, NJ},
        date={2016},
        ISBN={978-1-119-06195-3},
      review={\MR{3524748}},
}

\bib{etc}{inproceedings}{
      author={Botler, F{\'a}bio},
      author={Colucci, Lucas},
      author={Kohayakawa, Yoshiharu},
       title={The odd chromatic index of almost all graphs},
        date={2020},
   booktitle={Anais do {V} {E}ncontro de {T}eoria da {C}omputação},
   publisher={SBC},
     address={Porto Alegre, RS, Brasil},
       pages={49\ndash 52},
         url={https://sol.sbc.org.br/index.php/etc/article/view/11087},
}

\bib{botler21:_mod_chrom_index_rgs}{inproceedings}{
      author={Botler, F{\'a}bio},
      author={Colucci, Lucas},
      author={Kohayakawa, Yoshiharu},
       title={The $\!{}\bmod k$ chromatic index of random graphs},
        date={2021},
   booktitle={Extended abstracts---{EuroComb}~2021},
      editor={Ne{\v{s}}et{\v{r}}il, Jaroslav},
      editor={Perarnau, Guillem},
      editor={Ru{\'e}, Juanjo},
      editor={Serra, Oriol},
   publisher={Springer International Publishing},
     address={Cham},
       pages={726\ndash 731},
}

\bib{BoCoKo20}{article}{
      author={Botler, F{\'a}bio},
      author={Colucci, Lucas},
      author={Kohayakawa, Yoshiharu},
       title={The $\!{}\bmod k$ chromatic index of graphs is~{$O(k)$}},
        date={2023},
     journal={J. Graph Theory},
      volume={102},
      number={1},
       pages={197\ndash 200},
}

\bib{hasanvand22:_modul}{article}{
      author={Hasanvand, M.},
       title={Modulo factors with bounded degrees},
        date={2022},
     journal={arXiv e-prints},
      eprint={2205.09012},
}

\bib{haxell95ramsey}{article}{
      author={Haxell, P.~E.},
      author={Kohayakawa, Y.},
      author={\L{}uczak, T.},
       title={The induced size-{R}amsey number of cycles},
        date={1995},
        ISSN={0963-5483},
     journal={Combin. Probab. Comput.},
      volume={4},
      number={3},
       pages={217\ndash 239},
         url={https://doi.org/10.1017/S0963548300001619},
      review={\MR{1356576}},
}

\bib{Janson2000random}{book}{
      author={Janson, Svante},
      author={\L{}uczak, Tomasz},
      author={Ruci\'nski, Andrzej},
       title={Random graphs},
      series={Wiley-Interscience Series in Discrete Mathematics and
  Optimization},
   publisher={Wiley-Interscience, New York},
        date={2000},
        ISBN={0-471-17541-2},
         url={https://doi.org/10.1002/9781118032718},
      review={\MR{1782847}},
}

\bib{pyber1991covering}{inproceedings}{
      author={Pyber, L.},
       title={Covering the edges of a graph by {$\dots$}},
        date={1992},
   booktitle={Sets, graphs and numbers ({B}udapest, 1991)},
      series={Colloq. Math. Soc. J\'{a}nos Bolyai},
      volume={60},
   publisher={North-Holland, Amsterdam},
       pages={583\ndash 610},
      review={\MR{1218220}},
}

\bib{scott1997graph}{article}{
      author={Scott, A.~D.},
       title={On graph decompositions modulo {$k$}},
        date={1997},
        ISSN={0012-365X},
     journal={Discrete Math.},
      volume={175},
      number={1-3},
       pages={289\ndash 291},
         url={https://doi.org/10.1016/S0012-365X(96)00109-4},
      review={\MR{1475859}},
}

\end{biblist}
\end{bibdiv}
\endgroup %
\end{document}